%% file: ms.tex
\documentclass[letter,leqno]{article}

\usepackage[thinspace,amssymb,thinqspace,thinspace]{SIunits}
\usepackage{amsmath}
\usepackage{stmaryrd}
\usepackage{amssymb}
\usepackage{graphicx}
\usepackage{subfigure}
\usepackage{booktabs}
\usepackage{xspace}
\usepackage{cite}
\usepackage{soul}
\usepackage{multirow}

\usepackage{algorithm}
\usepackage[algo2e]{algorithm2e}
\usepackage[]{algpseudocode}

\usepackage[left=3cm,right=2cm,top=2cm,bottom=2cm]{geometry}
\usepackage[utf8]{inputenc}

\usepackage{tikz}

\usepackage[citecolor=black]{hyperref}

\include{tikz_defines}

\DeclareMathOperator*{\argmin}{arg\,min}
\DeclareMathOperator*{\argmax}{arg\,max}
\def\linefill{%
\leavevmode
\leaders\hrule\hskip\dimexpr\textwidth-2pt\mbox{}}

\newcommand{\optP}{\widetilde{P}}
\newcommand{\revisedP}{\bar{P}}
\newcommand{\filteredA}{\bar{A}}
\newcommand{\ddom}{{\cal D}}
\newcommand{\twiceFilteredA}{\widetilde{A}}
\newcommand{\filteredD}{\ensuremath{\bar{D}}\xspace}
\newcommand{\Ptent}{\ensuremath{{P}^{(t)}}\xspace}
\newcommand{\maxFilEval}{\bar{\lambda}_m}
\newcommand{\maxFilEvalEll}{\bar{\lambda}_{\ell,m}}
\newcommand{\maxEval}{{\lambda}_m}

\newcommand{\RsumD}{\widetilde{D}}

\newcommand{\maxRsumEval}{\widetilde{\lambda}_m}
\newcommand{\REMOVE}[1]{}

\newcommand{\rootstencil}{\textsf{Sprsfy}\xspace}
\newcommand{\spreadlumping}{\textsf{OffLmp}\xspace}
\newcommand{\dmod}{\textsf{1Norm}\xspace}
\newcommand{\econstraint}{\textsf{Cnstrnt}\xspace}
\newcommand{\codeComment}[1]{/\hskip -.05in / #1}

\definecolor{darkgreen}{rgb}{0,0.40,0}

\graphicspath{{fig/}{eps/}}

\begin{document}

\title{Smooth Aggregation for Difficult Stretched Mesh and Coefficient Variation Problems\thanks{
Technical report SAND2021-3298 O.
{This work was supported by the U.S.~Department of Energy, Office of Science, Office of Advanced Scientific
       Computing Research, Applied Mathematics program.
       Sandia National Laboratories is a multimission laboratory managed and operated by National Technology
       and Engineering Solutions of Sandia, LLC., a wholly owned subsidiary of Honeywell International, Inc.,
       for the U.S. Department of Energy's National Nuclear Security Administration under grant~DE-NA-0003525.
       This paper describes objective technical results and analysis.  Any subjective views or opinions that
       might be expressed in the paper do not necessarily represent the views of the U.S. Department of Energy
       or the United States Government.}}}
\author{
  Jonathan J. Hu\thanks{Sandia National Laboratories, P.O. Box 969, MS 9159, Livermore, CA 94661
    ({jhu@sandia.gov}, {rstumin@sandia.gov}).}
  \and
  Chris Siefert\thanks{Sandia National Laboratories, P.O. Box 5800, MS 1320, Albuquerque, NM 87185
    ({csiefer@sandia.gov}).}
  \and
  Raymond S. Tuminaro\footnotemark[2]
}

\maketitle

\begin{abstract}
\input{abstract}

\end{abstract}

\input{intro}

\input{algoDescription}

\input{results}

\input{conclusion}

\bibliographystyle{plain}

\input{ms.bbl}
\appendix

\input{appendix_randcube}
\input{appendix_stretchcube}

\end{document}

%% file: abstract.tex
Four adaptations of the smoothed aggregation algebraic multigrid (SA-AMG) method are 
proposed with an eye towards improving the convergence and robustness of the solver
in situations when the  discretization matrix contains many weak connections.
These weak connections can cause higher than expected levels of fill-in  within the 
coarse discretization matrices and can also give rise to  sub-optimal smoothing
within the prolongator smoothing phase. These smoothing drawbacks
are due to the relatively small size of some diagonal entries within the filtered 
matrix that one obtains after dropping the weak connections. 
The new algorithms consider modifications to the Jacobi-like step that defines the 
prolongator smoother, modifications to the filtered matrix, and also direct
modifications to the resulting  grid transfer operators. 
Numerical results are given illustrating 
the potential benefits of the proposed adaptations.

%% file: intro.tex
\section{Introduction} \label{sec: intro}

The smoothed aggregation algebraic multigrid (SA-AMG) algorithm  was originally proposed over twenty years ago as an effective and scalable
solution strategy for linear systems arising from discretized elliptic partial differential equations (PDEs)~\cite{VaMaBr96,VaBrMa01}. 
The basic form of the original algorithm has been employed without major mathematical modification within numerous applications to tackle a wide 
range of non-trivial problems.  While generally successful, convergence difficulties can arise for some complex applications, even for 
matrices coming from elliptic PDEs. 

Multigrid methods, including SA-AMG, are based on the idea that simple relaxation methods such as Jacobi generally smooth high frequency errors 
and that these smoothed errors can then be accurately represented and more efficiently reduced by relaxation iterations on a coarser grid 
representation of the linear system.  This paper considers modifications to the SA-AMG method, targeting some potentially vulnerable
components that tend to be more fragile when irregular or anisotropic coarsening is needed.  While simple relaxation methods
generally smooth errors, they do not necessarily do so in a uniform or isotropic fashion. That is, errors after relaxation may be 
much smoother in certain directions than in others. In these cases,  algebraic multigrid (AMG) coarsening 
must only occur in directions where errors are {\it algebraically smooth}
(i.e., directions where the relaxation method significantly damps some local error components).
These algebraically smooth errors are not necessarily geometrically smooth, but they do need to be well
represented within the range of the coarse grid interpolation.

Irregular coarsening occurs when the  strength-of-connection phase of the standard SA-AMG algorithm labels many nonzeros as {\it weak}. 
The presence of many weak connections accentuates some problematic facets of the 
hierarchy construction process. Specifically, 
SA-AMG constructs a graph based on only the strong matrix connections
(as determined by a strength-of-connection algorithm), which is then coarsened to generate coarse discretization operators. 
In general, sparser graphs (i.e., graphs with few connections) lead to less coarsening, larger coarse matrices,
denser coarse matrices, and more total AMG levels. The overall effect
is a potentially significant cost increase in both the 
setup and apply phases. Further, the presence of weak connections requires that a special filtered discretization
matrix be devised in order to define the grid transfers. This filtering step can be somewhat fragile, and if not done properly
can lead to sub-optimal grid transfer operators. All of this is exacerbated by limitations in the  standard strength-of-connection 
algorithms.  The original strength-of-connection idea dates back to the 1980s~\cite{BrMcRu84,RuSt85} and is 
motivated by M-matrix assumptions.  While computationally inexpensive, it is easy to construct examples
(see \S~\ref{sec:motivation}) where common strength criteria give a poor indication of 
the 
directions that correspond to algebraically smooth errors.  To avoid the adverse convergence effects of coarsening in 
directions that are not aligned with algebraically smooth errors, often a somewhat large threshold or cutoff value is chosen to 
encourage the labeling of many connections as weak, and thus coarsen slowly.
While there are some interesting alternatives to the classical coarsening 
approaches~\cite{Livne04,OBroeker_2003a,BrZi2005,BrBrMaMaMc06,BrBrKaLi15,BrannickF10,OlScTu10},
that might reduce some of these ill effects,
most of these alternatives are expensive and not fully robust. For this reason, these
alternatives have not been generally adopted and are not considered further in this paper. 

Four algorithmic adaptations are proposed to address potential SA-AMG deficiencies that tend to arise in the presence of many 
weak connections.  The new algorithms are algebraic and do not require any additional information or intervention from the application. 
Furthermore, these adaptations fit naturally into the existing SA-AMG setup workflow.  
The first idea considers an alternative diagonal approximation to the matrix inverse used within the Jacobi 
{\em prolongator smoothing} step.
Relatively small diagonal entries  can occur when formulating the filtered discretization operators. Unfortunately, these
small entries are highly problematic for the Jacobi step.
Effectively, we propose an alternative Jacobi-like step that uses the inverse 
of a diagonal matrix based on the 1-norm of individual rows of the filtered matrix.
The second algorithm addresses how the filtered matrix is defined. Normally, the filtered diagonal entries
are modified to reflect weak off-diagonal entries that are dropped from the original matrix. To prevent potentially small
diagonal entries, we propose that dropped entries are accounted for by modifying off-diagonal nonzeros in some circumstances.  
The third variant introduces a set of constraints that the prolongator must normally satisfy
(e.g., all entries lie between zero and one).  The algorithm attempts to reformulate any SA-AMG prolongator row that violates 
the constraints by finding a suitable nearby row.
The fourth and final algorithm introduces a second dropping or filtering stage
that is performed after the standard SA-AMG coarsening algorithm. That is, the second filtering stage
does not alter the coarsening process, but indirectly yields a sparser grid transfer operator
by effectively sparsifying the matrix used within the prolongator smoothing step. Here, the aim is
to reduce the cost of the solver
while maintaining SA-AMG's convergence properties.
In some cases, these four modifications do not
significantly alter cost or convergence behavior (when compared to traditional SA-AMG).  Under certain conditions, however,
modified SA-AMG performs much better than the traditional version, especially for problems where
irregular or anisotropic coarsening is needed.

The paper is structured as follows.  In \S\ref{sec: saamg} we give an overview of the SA-AMG method.  In
\S\S\ref{sec:1norm diag approx}--\ref{sec:root_node},
we present the new algorithmic variants for
how SA-AMG generates the grid transfer.
In \S\ref{sec:results_randcube} and \S\ref{sec:results_stretchcube} we present results for Poisson and
reaction-diffusion tests problem on a regular domain. In \S\ref{sec:results_sp10}, we explore the
effectiveness of the new algorithms on systems arising from a well-known reservoir model.  In \S\ref{sec:exawind}, we detail experiments
on systems arising in a low-Mach computational fluid dynamics wind turbine simulation.

\subsection{Smoothed Aggregation} \label{sec: saamg}

\input{background}

\input{dbar}

%% file: background.tex
An example multigrid V cycle iteration is given in Algorithm \ref{multigrid code} to solve
\begin{equation}
  A_1 u_1 = f_1. \label{eqn:Ax=b}
\end{equation}
\begin{algorithm}
\centering
\begin{tabbing} 
\hskip .6in \=  \hskip .3in \=  \hskip .3in \=  \hskip .3in \=  \kill  \\
  \> Vcycle($ A_\ell, f_\ell, u_\ell, \ell $) \{ \\
\>  \> if $\ell \ne N_{level} $ \{ \\
\>  \> \> $ u_\ell = \hat{S}_\ell (A_\ell, f_\ell, u_\ell )$ \\
\>  \> \> Vcycle($A_{\ell+1}, P_{\ell}^T (f_\ell - A_\ell u_\ell), 0, \ell\hskip -.04in + \hskip -.04in1 $) \\
\>  \> \> $ u_\ell = u_\ell + P_{\ell}  u_{\ell+1}$ \\
\>  \> \> $ u_\ell = \hat{S}_\ell (A_\ell, f_\ell, u_\ell )$\\
\>  \>  \} \\
\>  \> else  $ u_\ell = A_\ell^{-1} f_\ell $ \\
\> \}
\end{tabbing} 
\caption{\label{multigrid code}$N_{level}$ multigrid V-cycle to solve $A_\ell  u_\ell = f_\ell$.}
\end{algorithm}
To fully specify the algorithm, one must define relaxation procedures $\hat{S}_\ell()$, $\ell=1,\dots,N_{level}$
and the grid transfer operators $P_\ell$, $\ell=1,\dots,N_{level}-1$. Operator
$P_\ell$ interpolates solution updates from level $\ell\hskip -.04in+\hskip -.04in1$ to level $\ell$, while
operator $P^T_\ell$ restricts residuals from level $\ell$ to level $\ell\hskip -.04in+\hskip -.04in1$.
The coarse grid discretization operators $A_{\ell+1}$ ($\ell \ge 1$) are defined 
by 
\begin{equation} \label{Galerkin}
  A_{\ell+1} = P_{\ell}^T A_\ell P_{\ell}.
\end{equation}
Smoothed aggregation is a particular type of AMG method\cite{VaMaBr96,VaBrMa01,MaBrVa1998}
to determine the  $P_{\ell}$'s. Specifically, $P_{\ell}$ is given by 
\begin{equation} \label{eq:prolongator_smoothing}
    P_\ell = \left (I - \omega \filteredD_\ell^{-1} \filteredA_\ell \right) \Ptent_\ell,
  \end{equation}
where the term in parenthesis is referred to as the  {\em prolongator smoothing step },
$\omega =  \frac{4}{3 \maxFilEvalEll}$ is the {\em damping factor}, and $\Ptent_\ell$ 
is the {\em tentative prolongator}. Here, $\filteredA_\ell$ is a {\em filtered} form of $A_\ell$ with
diagonal $\filteredD_\ell$, and $\maxFilEvalEll$ is an approximation
to the maximum eigenvalue of $\filteredD_\ell^{-1} \filteredA_\ell$.
By filtered, we mean that some entries are dropped, as described shortly.
When all matrix connections are labeled as strong, $\filteredA_\ell =A_\ell $.

The tentative prolongator $\Ptent_\ell$ must accurately interpolate certain near null space (kernel) components of the discrete operator $A_\ell$. 
The prolongator smoothing step then improves the grid transfer operator by smoothing the basis functions associated 
with \Ptent. When many off-diagonals of $A_\ell$ are labeled as weak, $\filteredA_\ell \not=A_\ell $, and so a suitable $\filteredA_\ell$
must be defined. Before doing this, we first describe the coarsening process that defines $\Ptent_\ell$. Though smoothed aggregation can be 
applied to PDE systems, this paper focuses on scalar PDEs. In this case, $\Ptent_{\ell}$ is an $N_\ell \times N_{\ell+1}$ matrix with 
its sparsity pattern determined by a decomposition of the set of
$A_\ell$'s graph vertices into $N_{\ell+1}$ disjoint \textit{aggregates} $\mathcal{A}_\ell^i$, such that
\begin{equation} \label{eq:aggregates}
  \bigcup_{i=1}^{N_{\ell+1}} \mathcal{A}_\ell^i = \left\{ 1,...,N_\ell  \right\} \ , \ 
  \mathcal{A}_\ell^i \cap \mathcal{A}_\ell^j = \emptyset \ , \ i \neq j  \ .
\end{equation}
An ideal aggregate $\mathcal{A}_\ell^i$ is formed by grouping together a central or root vertex with its immediate neighbors (in the graph of $\filteredA_\ell$)
when both the central vertex and the neighbors are all unaggregated.  While ideal aggregates would only consist of a root vertex and its immediate neighbors, 
it is generally impossible to partition all of the vertices into ideal aggregates. Thus, some heuristics are needed to either create smaller aggregates
or enlarge the ideal aggregates. Each aggregate on level $\ell$ gives rise to one vertex on level $\ell+1$.

The nonzero values of $\Ptent_\ell$ are defined by partitioning the near null space of $A_\ell$ over the aggregates. In this paper and for scalar PDEs, the near null space
is simply the vector of all ones. This is related to the fact that constant functions lie in the null space of 
second order PDEs with the general form
$$
( a^{(1)}( {\bf x})~ u_x )_x  + ( a^{(2)}({\bf x})~ u_y )_y  + ( a^{(3)}({\bf x})~ u_z )_z  = f   \hskip 1in {\bf x} ~\in~ \Omega
$$
when only Neumann conditions are applied on the boundary of $\Omega$. Here,  the $a^{(k)}( {\bf x})$ are functions in $H^1 (\Omega)$. 
With this null space, $\Ptent_\ell$ is given by 
\begin{equation} \label{eq:aggregate identity}
  \left ( \Ptent_\ell \right )_{ij} = \left\{ 
\begin{array}{cl}
$1$ & \mbox{if} \ i \in \mathcal{A}_\ell^j \\
$0$ & $otherwise$  . \\
\end{array} \right.
\end{equation}
Near null space error components are not damped by conventional relaxation procedures and so this choice of $\Ptent_\ell$ guarantees
that these constants are well represented on coarse levels.
If the null space of $\filteredA_\ell$ is also defined by the space of constant vectors, then it is clear that the space of 
constants is contained not only in the range space of $\Ptent_\ell$ but also in the range space of $P_\ell$.

To complete the description of traditional SA-AMG for scalar PDEs, the filtered version $\filteredA_\ell$ of $A_\ell$ is given 
by 
\begin{equation}\label{eq:filteredA}
\left ( \filteredA_\ell \right )_{ij} = \left \{
\begin{array}{ll}
                           \left ( A_\ell \right )_{ij}  & \mbox{if~~} i \ne j \mbox{~and~}   j\in {\cal S}_{\ell,i}\\
                           \left ( A_\ell \right )_{ii} + \sum_{k \in {\cal W}_{\ell,i}} \left ( A_\ell \right )_{ik} & \mbox{if~~}  i = j \\
                           0       & \mbox{otherwise}
\end{array}
\right .
\end{equation}
where ${\cal S}_{\ell,i}$ and ${\cal W}_{\ell,i}$ are the sets of strong (weak) neighbors of the $i^{ih}$ vertex on level $\ell$.
With this definition, the sum of the entries within each row of $A_\ell$ is equal to the
sum of the entries within the corresponding $\filteredA_\ell$ row. Notice that if
the constant is within the null space of $A_\ell$ (i.e., the sum of entries within a row is zero),
it is also within the null space of $\filteredA_\ell$.
Finally, off-diagonals nonzeros are normally considered strong if and only if they satisfy 
\begin{equation}\label{eq:crit_sa}
| \left ( A_\ell \right )_{ij}| \geq \theta \sqrt{ \left ( A_\ell \right )_{ii} \left ( A_\ell \right )_{jj}},
\end{equation}
for some user-defined threshold $\theta\in[0,1]$.

%% file: dbar.tex
\subsection{Prolongator Smoothing and Small Diagonal Entries }\label{sec:motivation}

We now motivate \S\ref{sec:sa variations} by considering some possible prolongator smoothing shortcomings.
The prolongator smoothing step
corresponds to one iteration of the damped Jacobi method applied to the matrix equation 
\begin{equation} \label{eq:energy minimization}
\filteredA P = 0,
\end{equation}
starting with an initial guess of $P=\Ptent$ and $\filteredA$ given by \eqref{eq:filteredA}.
Here, we now drop the subscript $\ell$ to simplify the notation for the remainder of the paper.
Of course, $P = 0$ is a trivial solution
as is $P = C$ when  $\filteredA$'s null space is given by the space of constant vectors. Here, $C$ is a 
$n \times m$ dense matrix with
$c_{ik} = c_{q  k}$ for all $i$ and $q$ such that $ 1 \le i,q \le n $ and each $k$ such that $1 \le k \le m$.
While repeated Jacobi iterations 
converge to something uninteresting,
one Jacobi step extends the nonzero support of the interpolation basis functions by one in the graph of $\filteredA$.
When suitably damped, it also reduces the energy of these basis functions in the norm defined by 
$||p||_{\filteredA} = \sqrt{p^T \filteredA p}$ when $\filteredA$ is symmetric positive definite.  
Most AMG convergence theories rely on some bound for the energy of the interpolation basis functions. 
$\filteredD^{-1}$ can be viewed as an inexpensive approximation to
$\filteredA^{-1}$ that will reduce high frequencies in $\Ptent$'s basis functions when used
within a Jacobi step. 

Undamped Jacobi is guaranteed to converge when $\filteredA$ is strictly diagonally dominant, i.e., the magnitude of the diagonal
entry in every row is greater than the sum of the absolute value of all other nonzeros in that row. 
When $\filteredA$ is symmetric positive definite but not strictly diagonally dominant, damping may 
be needed to ensure that the spectral radius
of $I - \omega \filteredD^{-1} \filteredA $ is less than one.  For smoothed aggregation the damping parameter is 
typically chosen as $\omega = 4 / ( 3 \maxFilEval)$ and is based on a Chebyshev minimization principle~\cite{Br1997}.
For strictly diagonally dominant matrices, it is easy to see that $\maxFilEval < 2$ using the 
Gershgorin circle theorem. When the matrix is not strictly diagonally dominant, 
it is possible that $\maxFilEval \gg 2$ or even undefined if some $\filteredD_{ii} = 0$.
Unfortunately, $\maxFilEval $ might be large even if only one row violates the strictly diagonally
dominant condition.  Obviously, the Jacobi step has little effect when $\omega$ is small,
which occurs for large $\maxFilEval $. 
Some discretization schemes can produce diagonally dominant matrices\footnote{Diagonally dominant 
matrices are defined by replacing the {\it greater than the sum}  condition in the strictly diagonally 
dominant definition by a condition {\it greater than or equal to the sum}.}
where $\maxEval \le 2$.
However, matrices arising from finite element discretization are generally not diagonally dominant, 
though often $\maxEval$ is still not much larger than 2. Unfortunately, 
dropping might lead to a
$\maxFilEval$ that is much larger than $\maxEval$. In fact, a $\filteredD_{ii}$ might even 
become zero or negative.
For example, a nodal linear finite
element discretization of a Poisson operator on a mesh where all elements are 
$1 \times 1 \times h_z$ hexahedrons produces matrices with 27-point
interior stencils. After scaling by $36 h_z$, these
identical stencils are given by 
$32+64 h_z^2$ (diagonal entry), 
$-1 - 2 h_z^2$ (8 neighbors on cell corners),
$-4- 2 h_z^2$ (8 neighbors sharing a $x-z$ or $y-z$ face), 
$16 h_z^2 - 16$ (2 neighbors sharing a $z$ edge), 
$2 - 8 h_z^2$ (4 neighbors sharing a $x-y$ face), 
and $8-8 h_z^2$ (4 neighbors sharing a $x$ or $y$ edge).
When $h_z = \sqrt{7}/2$, these stencils take on the values 
$144 , -4.5 ,  -7.5 ,  12 , -12, ~\mbox{and}~ -6$.
For  $h_z > \sqrt{7}/2$, the positive off-diagonal entry becomes the largest in magnitude value. Thus, any magnitude-based
dropping criteria that employs a threshold that happens to retain the largest entry but drops all others will 
result in a negative
$\filteredA_{ii}$. This follows from the fact that $\filteredA$'s rows all sum to zero.

Further, small diagonal entries can be even more problematic for non-symmetric systems. 
While there may not be a damping parameter to guarantee Jacobi convergence
for non-symmetric matrices, it is still effective for matrices that are {\it close} to
symmetric positive definite.  In a wind simulation that will be
described in \S~\ref{sec:exawind}, 
the discretization matrix includes a symmetric sub-block of only slightly smaller dimension than the
entire matrix. However, a few of the non-symmetric rows/columns have small diagonal entries, 
so that the non-symmetric part of $\filteredA$ is emphasized in the matrix $\filteredD^{-1} \filteredA$.
We have observed that this causes 
severe eigenvalue convergence problems and very poor $\maxFilEval$ estimates 
coming from the typical power method used to provide these estimates.

%% file: algoDescription.tex
\section{Smoothed Aggregation Variations} \label{sec:sa variations}
We now outline four algorithm variations to address some of the difficulties associated with weak connections
and small aggregates.

\subsection{A 1-norm diagonal approximation}\label{sec:1norm diag approx}

As noted, the $\filteredD^{-1}$ within the Jacobi iteration can be viewed as a diagonal approximation to 
$\filteredA^{-1}$.  There are, however, alternatives, and one natural possibility is to replace $\filteredD$ by
the diagonal matrix $\RsumD$ where 
$\RsumD_{ii}$ is the sum of the absolute values of nonzeros
in the $i^{th}$ row of $\filteredA$.
That is,
\begin{equation} \label{eq:absRowSum}
\RsumD_{ii} = \sum_{j} | \filteredA_{ij} | ~~,
\end{equation}
which is the $1-$norm of this $i^{th}$ row.
One could argue that this $1-$norm 
choice for $\RsumD_{ii}$ better captures the scaling of the entire row and that an iterative procedure based on 
this will be less sensitive to the diagonal dominance properties of the matrix.  This approximation 
is not new, and is related to a variant of the well-known pressure-correction
algorithm, SIMPLE, that is referred to as SIMPLEC~\cite{doi:10.1080/01495728408961817}. As we will see, this relatively straight-forward change can have a 
significant convergence effect.
Notice that 
the entries of this $\RsumD$ matrix are larger (assuming more than one nonzero per row) than those of $\filteredD$, and so the entries of $\RsumD^{-1} \filteredA$ are smaller in magnitude than those of $\filteredD^{-1} 
\filteredA$.  Overall, however, it is the entries of $\frac{4}{3 \maxFilEval} \filteredD^{-1} \filteredA$ and 
$\frac{4}{3 \maxRsumEval} \RsumD^{-1} \filteredA$ that appear within the Jacobi iteration, and these will be
comparable when $\filteredA$ is diagonally dominant, as $\maxRsumEval $ will be smaller than $\maxFilEval $. 
For example, when every row of $\filteredA$ is given by the stencil
$$
\begin{matrix}
      & -1  &      \cr
    -1  &  ~~4  &   -1   \cr
     & -1  &    
\end{matrix} ~~,
$$
then it is easy to see that 
$\frac{1}{\maxFilEval} \filteredD^{-1} \filteredA = \frac{1}{\maxRsumEval} \RsumD^{-1} \filteredA$.
This case corresponds to a constant coefficient periodic Poisson problem where 
all diagonal entries of $\filteredD$ are 4, all diagonal entries of $\RsumD$ are 8, 
$\maxFilEval = 2$, and $\maxRsumEval = 1$.  Clearly, this exact equality will hold whenever 
$\filteredA$ is defined by any circulant matrix.
Thus, the proposed $1-$norm diagonal modification will exactly reproduce the behavior
of the traditional smoothed aggregation method in this case. We normally expect to see similar SA-AMG convergence
behavior with either $\filteredD$ or  $\RsumD$ when $\filteredA$ is diagonally dominant. 
For more general problems, however, we will see that the behavior can be significantly different.

One interesting consequence of this $\RsumD$ definition is that the magnitude of $\maxRsumEval $,
the maximum eigenvalue associated with $\RsumD^{-1} \filteredA$, is always bounded by 1,
again as a consequence of the Gershgorin circle theorem. This implies that one can 
omit the eigenvalue calculation (which is problematic for the wind
simulation in \S~\ref{sec:exawind}) and consider using
$1$ as an estimate for $\maxRsumEval $.

\subsection{A Safe Guarded Diagonal Approximation}\label{sec:safeguard}

While this $\RsumD_{ii}$ choice works well in practice, we propose one additional modification that is
motivated by cases when the diagonal entry within a row of $\filteredA$ is much larger than the 
sum of the magnitudes of the row's off-diagonal entries.  
While this cannot happen for rows where the sum of
the entries is zero,
it 
might occur for matrix rows associated with or near Dirichlet boundary conditions.
To better understand this large diagonal scenario, consider the extreme case
when $\filteredA$ is in fact a diagonal matrix. Then, $\RsumD^{-1} \filteredA  = I $, $\maxRsumEval$ is 1, 
and the prolongator smoothing step \ref{eq:prolongator_smoothing} becomes
$$
P = \left (I - \frac{4}{3\lambda_m}\RsumD^{-1} \filteredA   \right) = \left (I - \frac{4}{3} I \right) \Ptent  =  -\frac{1}{3} \Ptent .
$$
If instead the damping parameter had been defined as $\omega = 1 / ( \maxRsumEval)$, which is
not an unreasonable choice, then $P$ would be identically zero.
Essentially, the Jacobi iteration converges too rapidly to the $P = 0$ solution when the diagonal is
much larger than the sum of the magnitudes of the row's off-diagonals. Instead of only reducing high frequencies,
the whole spectrum is reduced, which is undesirable. If only a subset
of rows have a very large diagonal, it would be better
to use the corresponding row of $\Ptent$ or to somehow limit the damping. While this could be
accomplished by setting $\omega$ to a suitably small value, this would have the unintended 
consequence of limiting the smoothing effect for \textit{all} rows, even those that are not so diagonally
dominant. Instead, we prefer modifying or boosting the value of $\RsumD_{ii}$  for any row
where the modified Jacobi step significantly reduces low frequencies. This effectively corresponds
to augmenting the definition of $\RsumD_{ii}$ with a safeguard. To do this, we must have a
criterion to detect the low frequency reduction.

To understand a possible remedy for low frequency reduction, 
consider the effect of the prolongator smoother step on just one row of the matrix. Specifically, assume that
vertex $i$ defines the center or root point of the $k^{th}$ aggregate. That is, $i \in {\cal A}_k $  and 
$j \in {\cal A}_k $ for any $j$ such that $\filteredA_{ij} \neq 0$ . If we additionally take $\maxRsumEval = 1$, 
then the $i^{th}$ row of the resulting prolongator is given by 
$$
P_{ij} = 1 - \frac{4}{3} \frac{s_i}{\RsumD_{ii}},
$$
where $s_i$ denotes the sum  of $\filteredA$'s entries in the $i^{th}$ row.
Note that if $s_i=0$ and $\RsumD_{ii} \neq 0$, then $P_{ij} = 1$. One can argue that it is natural for a prolongator 
basis function to have a value of $1$ at the aggregate's central node  and to then decay smoothly toward
zero.  This is analogous to the use of injection to interpolate points that are co-located on the fine
and coarse grids. Here, a coarse unknown associated with an aggregate is viewed as being co-located 
with the root node. As the basis function should be largest at the root node, it should not have a negative
value or be too small at this root vertex. 
In other words, we seek to enforce the condition that
\[
P_{ij} \ge \sigma,
\]
where $i$ is again the root node of the $j^{th}$ aggregate
and $\sigma$ is the minimum acceptable value for the prolongator basis function at the root vertex. 
This implies that $\RsumD_{ii}$ must be chosen
such that 
\[
1 - \frac{4}{3} \frac{s_i}{\RsumD_{ii}}  \ge \sigma,
\]
or 
\begin{equation} \label{eq:not too small}
  \frac{4 s_i}{3 (1 - \sigma)} \le \RsumD_{ii} .
\end{equation}
Notice that when $s_i = 0$ and $\sigma < 1$, any non-negative $\RsumD_{ii}$ will satisfy this condition. 
Thus, \eqref{eq:absRowSum} can be used to define $\RsumD_{ii}$.
When $s_i$ is not zero, an initial value computed via \eqref{eq:absRowSum} is checked to see
that it satisfies \eqref{eq:not too small}. If this second condition is not satisfied, \eqref{eq:absRowSum} 
is discarded and instead the smallest  $\RsumD_{ii}$ satisfying \eqref{eq:not too small} is used to 
define $\RsumD_{ii}$. 
While the precise choice of $\sigma$ is not obvious, we have found experimentally that 
$\sigma = 1/3$ works well.
For this choice of $\sigma$, \eqref{eq:not too small} becomes
\begin{equation} \label{eq:not 2 small}
  2 s_i \le \RsumD_{ii}
\end{equation}
Thus, we first compute $\RsumD_{ii}$ using \eqref{eq:absRowSum}. For $\RsumD_{ii}$ that are identically zero
(due to an entire zero row of $\filteredA$) we set $\RsumD_{ii} = 1$,  and for those remaining rows violating \eqref{eq:not 2 small}
we set $\RsumD_{ii} = 2 s_i$. 
While the arguments were motivated by considering $i$ to be a root node of an aggregate, we apply this
criteria for all rows in the matrix, regardless as to whether or not the node is a root node. At non-root nodes, the
associated $P_{ij}$ can be smaller than $\sigma$ due to the fact $s_i$ does not reflect the $\filteredA \Ptent$
product.

\subsection{Alternative lumping strategies} \label{sec:alt lumping}
As already noted, 
near null space vectors should be accurately represented in the range space of the interpolation operator. 
If the null space of $A$ 
coincides with that of $\filteredA$ and is also contained within the range space of $\Ptent$, 
then it
will additionally be contained within the range space of the smoothed prolongator.
To force the
two null spaces to coincide, $\filteredA$'s entries must be modified to account for dropped nonzeros.
This is traditionally accomplished by only changing the diagonal.
As discussed, any dropped off-diagonal $A_{ij}$ is simply added
to  $A_{ii}$ to define $\filteredA_{ii}$. Unfortunately, this simple process can drastically alter the properties 
of the resulting filtered matrix. In extreme cases, it is possible that $A_{ii}$ and $\filteredA_{ii}$ 
have opposite signs.  
More generally, the diagonal dominance properties of a row might change.  That is, 
we could have 
$$
\ddom(\filteredA_i) \gg \ddom(A_i) 
~~\mbox{where}~~
\ddom(A_i) =  \frac{\sum_{j\ne i} |A_{ij}|}{|A_{ii}|} .
$$
Here, $A_i$ refers to the $i^{th}$ row of $A$ and 
$\ddom()$ measures the magnitude of off-diagonal entries relative to the diagonal. 
Notice that $\ddom(A_i) = 1$ when all off-diagonals are negative and the sum of $A_i$'s entries are zero.
In this case, any reasonable lumping strategy will result in $\ddom(\filteredA_i) = 1$.
However, this will not be the case when there are both positive and negative off-diagonals or when $A_i$'s entries
do not sum to zero (e.g., at a boundary).
For the most part, small values of $\ddom()$ are preferred by a damped Jacobi iterative method 
(i.e., the prolongator smoothing step) as 
its convergence rate is generally more rapid for matrices with small off-diagonal entries relative to 
the diagonal entries. This suggests that it might be effective to consider a
strategy that restricts or limits the lumping of terms to the diagonal to 
avoid significant growth in the resulting $\ddom(\filteredA_i)$. There are many possible such strategies.
We now describe a scheme that is primarily oriented toward scalar diffusion-type PDEs, 
and so may not be appropriate for other operators. 

For scalar diffusion-like PDEs, one can argue that positive off-diagonal entries are somewhat irregular. 
To see this, consider the simplified prolongator smoother step $ u = (I - \filteredD^{-1} \filteredA$)v .
When all off-diagonals are negative and when the constant
is in the null space of $\filteredA$, then each $u_i$ is just a weighted average of the $v_i$'s within 
its immediate neighborhood.
In this way, the relaxation process mimics a diffusion process associated with a heat equation. 
When some off-diagonal entries are instead positive, then some weights will be negative, which no longer 
resembles a diffusion process.  This can be made more rigorous by consider the relationship 
between Jacobi iterations and time marching for ordinary differential equations~\cite{OlScTu10}.
Thus, a possible lumping algorithm might make decisions based on the sign of matrix entries. 
Our overall lumping strategy considers modifying or perturbing the retained positive off-diagonals, 
the diagonal, or the retained negative off-diagonals {\it in this order of preference}
when the perturbation is negative.  The general aim is to 
enforce $\ddom(\filteredA_i) \le \tau \ddom(A_i)$ where $\tau$ is a user-supplied growth factor. 

\begin{algorithm} \label{algo:spreadLumping}
\begin{tabbing} 
\hskip .6in \=  \hskip .3in \=  \hskip .3in \=  \hskip .3in \=  \kill  \\ \linefill \\
$\filteredA_{i} = $\textsf{Lump\_AvoidSmallDiag}( $A_i, {\cal R}_i, \tau$)\\[2pt]
~~Input: \\
~~~~~${A}_i$      \> $i^{th}$ row of matrix with entries to be dropped \\
~~~~~${\cal R}_i$ \> set of column indices in $i^{th}$ row to be removed\\
~~~~~$\tau$       \> tolerance indicating that $\ddom(\filteredA_i)$ should not exceed $\tau \ddom(A_i) $ \\
~~Output: \\
~~~~~$\filteredA_{i} $ \> matrix row where $\filteredA_{ij} = 0 ~\mbox{for}~ j \in {\cal R}_i$ and 
$\filteredA_i v = A v $ where $v$ is a constant vector 
 \\[-4pt] \linefill \\[5pt]
~1. Let $r_i    \leftarrow \sum_{k \in {\cal R}_i  } A_{ik} $ \\[2pt]
~2. {\bf if} $r_i > 0 ~{\bf then}~ \filteredA_{ii} \leftarrow  A_{ii} + r_i$  ~~~\codeComment{decreases $\ddom(\filteredA_i)$}\\
~3. {\bf else \{} \\
~4. \> Let ${\cal K}_i^+  \leftarrow \{ k ~|~~ A_{ik} > 0 ~\land~ k \neq i ~\land~ k \notin {\cal R}_i \} $ \\[2pt]
~5. \> Let ${\cal K}_i^-  \leftarrow \{ k ~|~~ A_{ik} < 0 ~\land~ k \neq i ~\land~ k \notin {\cal R}_i \} $ \\[2pt]
~6. \> Let $\kappa_i^+  \leftarrow \sum_{k \in {\cal K}_i^+} A_{ik} $;~~~$\kappa_i^-  \leftarrow \sum_{k \in {\cal K}_i^-} A_{ik}$;\\[2pt]
~7. \>  {\bf if} $|r_i| \le \kappa_i^+$ {\bf then}~ ${\filteredA}_{ij} \leftarrow {A}_{ij} (1 + \delta_{i})$ for $j \in {\cal K}_i^+$ where $ \delta_i \leftarrow r_i / \kappa^+ $\\
~8. \>  {\bf else \{} \\
~9. \>  \>  ${\filteredA}_{ij} \leftarrow 0$ for $j \in {\cal K}_i^+$ \codeComment{zero out the ${\cal K}_i^+$ by distributing a} \\[3pt]
10. \>  \>  $\hat{r}_i \leftarrow r_i + \kappa_i^+  $ ~~~~~~~~\codeComment{portion of $r_i$ ($= \kappa_i^+$) to them} \\
11. \> \>          {\bf if} ${\cal K}_i^- == \emptyset$ {\bf then} redistribute to ${\cal K}^+$ if possible or if not \\
12. \>\> ~~~~~~~~~~~~~~~~~~~~~~~~possible do not modify row $i$ and return \\
13. \> \>          {\bf else \{} \\
14. \> \>  \> find largest positive $r_i^* < \min(d_{ii},|\hat{r}_i|)$ such that   $\ddom(\filteredA_{i}) \le \tau \ddom(A_i) $\\
15. \> \>  \>          define  $\filteredA_{i}$ such that its only nonzero values are \\
16. \> \>  \>          ~~~~~~~~ $\filteredA_{ii} \leftarrow {A}_{ii} - r_i^* $\\
17. \> \>  \>          ~~~~~~~~ $\filteredA_{ij} \leftarrow {A}_{ij} (1 + \delta_{i})$ for $i \in {\cal K}_i^-$  \\
18. \> \>  \>          ~~ where $ \delta_i \leftarrow (\hat{r}_i + r_i^* )/ \kappa^- $ \\
19. \> \>          {\bf \}} \\
20. \> {\bf \}}  \\
21. {\bf \}}  \\ 
\end{tabbing}
\caption{Alternative lumping procedure for prolongator smoothing step}
\end{algorithm}
A detailed algorithm description is given in Figure~\ref{algo:spreadLumping}. The algorithm is 
supplied a set of indices ${\cal R}_i$ denoting the nonzero columns that should be {\it removed} 
from the $i^{th}$ row. The sum of these entries $r_i$ must be then distributed to the {\it kept} entries of $\filteredA_i$.
If this sum is positive, then only the diagonal is modified in \textsf{Line 2}
as this lowers $\ddom(\filteredA_i)$.  If instead $r_i$ is negative, then  more care is necessary. 
We first split the set of kept 
indices into two subsets ${\cal K}^+$ and ${\cal K}^-$ corresponding to 
entries that have positive $A_{ij}$ values or negative $A_{ij}$ values, respectively.
The sum of the nonzero values associated with these two sets is denoted by ${\kappa_i}^+$ and ${\kappa_i}^-$,
respectively. \textsf{Line 7} corresponds to the case when all of the lumping can be distributed
to the positive kept entries without creating any new negative entries. If this is impossible, 
we distribute a portion of $r_i$ equal to $-\kappa_i^+$ to the ${\cal K}^+$, effectively zeroing them out.
That is, we prefer not creating new negative entries as this might fundamentally change the equation's character.
For the remaining $\hat{r}_i = r_i + \kappa_i^+$,
we seek in \textsf{Line 14} the largest magnitude perturbation,
$|r^*|$, to the diagonal 
that does not violate the $\ddom()$ growth restriction.
The remaining $\hat{r}_i - r_i^*$ is then distributed proportionally to the negative kept off-diagonals. 
In many cases,
$r_i^*  = \hat{r}_i$, so $\hat{r}_i$ is lumped entirely to the diagonal. When this is not true, $r_i^*$ is generally given by
$$
r_i^* = \frac{\hat{r}_i + \kappa_i^- + \tau \ddom(A_i) A_{ii} }{1-\tau \ddom(A_i)}  .
$$
This is obtained by some algebraic manipulations after first setting $\ddom(\filteredA_i) = \tau \ddom(A_i)$ and 
recognizing that the sum of the absolute values of $\filteredA_i$'s off-diagonals is $(-\kappa_i^-) (1 + \delta)$ while its
diagonal is $A_{ii} + r_i^*$.

However, safe-guards must be added for situations where there is no suitable value of $r_i^*$
satisfying the $\ddom()$ growth restriction. This might occur if the set ${\cal K}^-$ is empty. 
If it is instead possible to satisfy the $\ddom()$ growth restriction by further lumping to ${\cal K}^+$,
then this is done even though these off-diagonals now become negative. Otherwise, if there are no kept 
off-diagonal entries in the row, 
we skip the perturbation entirely, no longer preserving the row sums in these problematic
rows.

\subsection{Prolongator Constraints} \label{sec:constraints}
From a geometric multigrid perspective, one can argue that the entries of $P$ should lie 
between $0$ and $1$ inclusive.  For example, consider the vector $u = P e^{(j)}$ where $e^{(j)} $ is a 
canonical basis vector with only one nonzero entry, the $j^{th}$ element, that is set to one. It is clear
that all entries of $u$ should be positive and not greater than one for any 
sort of geometric interpolation scheme (as opposed to an extrapolation scheme).  For smoothed aggregation,
this connection is a little less apparent. Obviously, the tentative prolongator $\Ptent$
satisfies these constraints when the null space of $A$ is given by the vector of all ones. 
As the objective of the prolongator smoother step is to produce low energy grid transfer basis functions,
it can also be argued that these smoothed basis functions should decay smoothly from the peak
value to zero. Thus, these basis functions should not be negative anywhere. 
Though somewhat less obvious, the peak will typically be either one or less than one.
Specifically, smoothed prolongator basis functions  and tentative prolongator basis functions
will coincide for vertices that are not part of the aggregate boundary. This is a property of the 
null space of $\filteredA$ being a constant function in the case of
the Laplace operator and was discussed for ideal aggregates in \S~\ref{sec:safeguard}.
The remaining nonzero basis function entries will typically
be less than one due to the energy minimization and smooth decay properties just mentioned.
When  instead prolongator values  do not lie between $0$ and $1$ inclusive, it is often
an indicator that the smoothing of some basis functions is sub-optimal.
In these cases, one can consider enforcing a condition that all $\revisedP_{ij}$ lie 
between $0$ and $1$. Specifically, one can construct a minimization problem for the $i^{th}$
row of $\revisedP$ 
$$
\left\{
  \begin{aligned}
    &\revisedP_i = \argmin_{\optP_i} \; || \optP_i - P_i ||_2 \\
    &\;\,~~~~~~\text{\small subject to}\hspace{1em} 
\optP_{ij} = 0 ~\mbox{if}~ P_{ij} = 0,~~~
\optP_{ij} \ge 0,~~~ \optP_{ij} \le 1,~~~
\optP_i v = P_i v 
  \end{aligned}
\right .
$$
where $v$ is the vector of all ones. 
That is, find a new prolongator row that is closest to $P_i$ and satisfies the two bound constraints, has
a sparsity pattern not extending beyond $P_i$'s pattern, and where the sum of $P_i$'s entries and $\revisedP_i$'s entries  are identical.
If the row sum of $P_i$ is negative, then there is no feasible solution to this minimization problem.
If the $P_i$ row sum is zero, then the only solution is that $\revisedP_i$ is identically
zero. There is also no feasible solution when the row sum is greater than the number
of nonzeros in the sparsity pattern of $P_i$.  When the minimization problem has a
feasible solution, it can be obtained by the procedure described in
Algorithm~\ref{alg:post process}.  The algorithm itself is relatively
straightforward.  On each pass of the while loop, we take the worst
constraint violators on both ends (both negative entries and entries which are
greater than one), pin them to their constraint values, and then split the combined change in values
among all remaining non-violating entries.  This will then be repeated
until no constraint violating entries remain.
When the constraints cannot be satisfied within the $i^{th}$ row, then we simple take $\revisedP_{i} = \Ptent_{i}$.
\begin{algorithm}
    \begin{algorithmic}
\caption{$\revisedP_i$=\textsf{Constrain\_One\_P\_Row}$(P_i)$ \label{alg:post process}}
\State Let $\revisedP_{i} = P_{i} $
\State Let $w = \{ k ~|~~ \revisedP_{ik} \neq 0  \} $ \;
\State \While{$(\min_k \revisedP_{ik} < 0 )~ \land ~(\max_k \revisedP_{ik} > 1 )$}{
\State  $\hat{k} \gets \argmin_{k} \revisedP_{ik} $\;
\State  $\tilde{k} \gets \argmax_{k} \revisedP_{ik}$\;
\State  $\delta \gets 0$\;
\State  \lIf{$\revisedP_{i\hat{k}}   ~<~ 0 $} { $\delta \gets \delta +\revisedP_{i\hat{k}  }$; ~~$\revisedP_{i\hat{k}  } \gets 0$; ~~$w \gets w \setminus \hat{k}$}\;
\State  \lIf{$\revisedP_{i\tilde{k}} ~>~ 1 $} { $\delta \gets \delta +\revisedP_{i\tilde{k}} - 1$; ~~$\revisedP_{i\tilde{k}} \gets 1$; ~~$w \gets w \setminus \tilde{k}$}\;
\State  $\revisedP_{ik} \gets \revisedP_{ik}  + \delta/|w| $ ~for~ $k \in w$\;
}
    \end{algorithmic}
  \end{algorithm}

\subsection{Further sparsification}\label{sec:root_node}
\begin{figure}[ht!]
  \centering
  \includegraphics[scale=0.3]{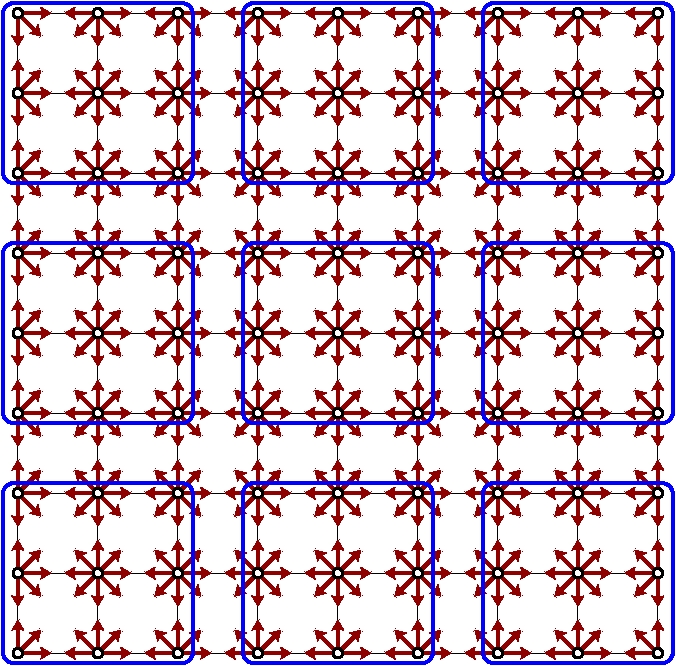} ~~~~~~~~~~~~
  \includegraphics[scale=0.3]{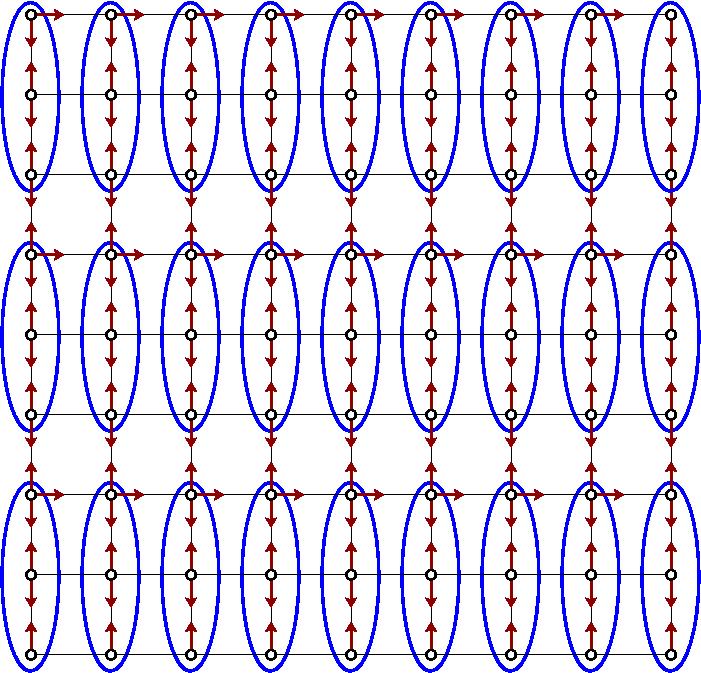}
  \caption{Two aggregation examples (aggregates indicated by blue enclosures) on a 2D structured mesh. Each red arrow indicate a strong $\filteredA_{ij}$ connection for vertex $i$ at the arrow base and vertex $j$ is the closest node in the direction pointed to by the arrow.}
  \label{fig:hotdog example}
\end{figure}

To understand a possible further sparsification of the smoothed aggregation prolongator operator, 
we first review the aggregation process within smoothed aggregation.
As noted, smoothed aggregation applies an algorithm to the graph of $\filteredA$ to construct aggregates
${\cal A}_j$ such that each fine mesh vertex belongs to only one aggregate. The basic idea is that a root node
is first chosen and then an initial aggregate is defined as the root node and all of its strong neighbors. Each root
node is chosen among vertices that have not yet been aggregated and are not adjacent (via strong connections) to 
any existing already
aggregated vertex.  This aggregation procedure is repeated until it is no longer possible to find such a root
node as all unassigned vertices are adjacent to assigned vertices. At this juncture, some heuristics are needed
to assign these remaining unassigned vertices by either creating new aggregates or enlarging existing aggregates.
This implies that most aggregates are 
composed of a central root node and its strong neighbors, thus the shape of the aggregates is primarily governed
by the strong neighbors of the root node. A typical aggregate will have a diameter of length 3. Figure~\ref{fig:hotdog example}
illustrates two aggregate scenarios on a regular mesh. In the leftmost image, all connections are strong 
and the corresponding aggregates happen to be perfect squares. In the rightmost image, 
only the vertical connections are strong with the exception of a few horizontal edges that in 
this contrived example never coincide with edges emanating from a root node. Once again, the aggregates are 
perfect, consisting of 3 points aligned in the vertical direction. That is, the aggregate shapes are
determined by the root nodes, which only have strong vertical connections. The main issue is that the horizontal
connections shown in the rightmost image will lead to nonzero fill-in within the coarse level discretization 
matrix due to the prolongator smoothing step.
Specifically, additional nonzeros connections arise between coarse vertices associated with non-neighboring aggregates (distance two aggregates) in the horizontal direction.
In particular, it is easy to show that for nontrivial $\filteredA$ matrices a nonzero $(i,j)$ entry 
occurs in the coarse matrix discretization whenever $(\Ptent_{.,i})^{T} \filteredA^T A \filteredA (\Ptent)_{.,j}$ is 
nonzero where $\Ptent_{.,j}$ refers to the $j^{th}$ column of $\Ptent$. This will certainly occur if there
is a distance 3 path in the filtered graph (along the red arrows in Figure~\ref{fig:hotdog example}). 
More generally, fill-in between distant aggregates can occur due to a conflict between the
characterization of strong connections between the root node and the non-root nodes. Specifically, the root node indicates that connections to certain neighboring aggregates
are weak
while some member of the root node's aggregate
has a strong connection to this very same aggregate. 
To reduce fill-in, we can look for conflicts and re-label some conflicting connections.

As our conflict characterization is based on aggregate choices, it is most practical to only
consider non-root strong connections for re-labeling 
as these do not alter the definition of the already-chosen aggregates.
This re-labeling corresponds to a further sparsification of $\filteredA$ that is performed immediately
preceding the prolongator smoother step. This sparsification occurs aggregate-by-aggregate. First, 
the weak connections of the root node are examined to determine the neighboring aggregates (in the graph 
of $A$) associated with these weak connections. If these weak-neighbor aggregates have no strong connections
to the root node, then these aggregates are put into a set {\sf CandidatesForPruning}.
Second, we examine all the non-root vertices in the aggregate looking for
strong connections to any vertex within each {\sf CandidatesForPruning} aggregate. If there is just one strong connection to a particular aggregate in {\sf CandidatesForPruning}, this strong connection is re-labeled as weak and dropped from the $\filteredA$ that will be used in the prolongator smoothing step. The diagonal entry is modified
to reflect the dropped entry following the usual method.
Thus, strong connections to each aggregate in {\sf CandidatesForPruning} are retained if there
are multiple strong connections but dropped if there is only one strong connection.
This leads to a further sparsification of $\filteredA$ that may now be non-symmetric even if 
$\filteredA$ is symmetric. 
To remedy this, all $(j,i)$ entries are dropped if the associated $(i,j)$ entry was dropped in the re-labeling phase.
This restores the symmetry in the resulting filtered
matrix, denoted as $\twiceFilteredA$, that is then used in the prolongator smoothing step. As the 
final prolongator is now sparser, we can expect that the amount of fill-in will be reduced on coarse
level matrices. However, this may also cause the convergence rate to be somewhat slower.

%% file: results.tex
\section{Numerical Results} \label{sec:results}
Four sample problems are presented to examine the behavior of the four
proposed SA-AMG algorithm variations. The first two correspond to
fairly academic cubes with perturbations or stretched mesh
spacing. The other two are more realistic. One is the SPE10 benchmark
problem~\cite{spe10} from 
the Society of Petroleum Engineers, which has highly heterogeneous material jumps. The second comes from a wind turbine 
simulation where the underlying mesh has problematic aspect ratios within some parts of the domain. 

\subsection{Randomly Perturbed Cube}\label{sec:results_randcube}
The first test problem considers a Poisson 
equation
\begin{equation} \label{LAPLACE_EQN}
\begin{aligned} 
-\nabla^2 u(x) &= f(x)   \hskip .4in  &x \mbox{~in~} \Omega, \\
u(x) &= u_D(x)                        &x \mbox{~on~} \partial \Omega,
\end{aligned}
\end{equation}
where the domain is defined by $\Omega=[0,1]\times[0,1]\times[0,100]$, the forcing function
$f(x)$ is identically zero, 
and $u_D(x) = 1 + x_1 + x_2 + x_3 + x_1 x_2+ x_1 x_3+ x_2 x_3 + x_1 x_2 x_3$. 
We discretize \eqref{LAPLACE_EQN} using linear hexahedral finite elements on
a tensor product mesh using 60 elements in each coordinate direction. 
The element sizes are given by a perturbation of a uniform spacing. Specifically, 
along each coordinate direction mesh points are randomly perturbed by up to $20\%$ 
from the uniform spacing location (details given in Appendix~\ref{appendix:randcube}).  
This leads to a mesh where the spacing in the first two coordinate directions is 
generally finer than that in the third dimension.
For our experiments, we generated 50 different test meshes using
different random seeds.

The smoothed aggregation multigrid solver in MueLu \cite{MueLu} is used as a preconditioner
for the conjugate gradient (CG) method to solve the 50 linear systems. It is well-known
that linear hexahedral elements on highly stretched meshes give rise to 
matrix coefficients that are problematic from a strength-of-connection perspective.
That is, the magnitude of the matrix entries are not well correlated with
mesh stretching. For this reason, an alternative matrix termed a distance Laplacian, 
$L^{(d)}$, is used in this experiment for the aggregation/coarsening phase of the algorithm. Specifically,
strong connections satisfy $ |L^{(d)}_{ij}| \geq \theta \sqrt{L^{(d)}_{ii} L^{(d)}_{jj}} $.
The matrix $L^{(d)}$ has the same nonzero pattern as the discretization matrix $A$. The $L^{(d)}_{ij}$
off-diagonal values, however, are defined as the reciprocal of the negative distance between the 
$i^{th}$ and $j^{th}$ coordinate, requiring coordinates be supplied to the solver. 
The diagonal is then chosen so that sum of all entries within each
row is identically zero. 
Here, $\theta$ is taken as $0.025$, which was determined
experimentally to produce desirable aggregates that are primarily oriented along the first two coordinate
directions. The smoothed 
aggregation hierarchy was generated so that the resulting discretization
matrix is coarsened until only 1,000 or fewer unknowns remain, at which point a
direct solver is applied.  Two sweeps of Chebyshev pre- and post-smoothing are applied on all 
other levels. The Chebyshev eigenvalue interval is given by $[\lambda^*/10,\lambda^*]$ where
$\lambda^*$ is an estimate of the maximum eigenvalue of $A$ obtained by $10$ sweeps
of the power method.  The conjugate gradient
iteration is terminated when the residual is reduced by a factor of $10^{-10}$.

We now consider all 16 combinations of enabling and disabling the four adaptations: 
the \dmod diagonal modification            of \S\S~\ref{sec:1norm diag approx}--\ref{sec:safeguard},
the \spreadlumping lumping modification    of \S~\ref{sec:alt lumping},
the \econstraint   constraint enforcement  of \S~\ref{sec:constraints}, and
the \rootstencil   sparsification          of \S~\ref{sec:root_node}.
Disabling all four options corresponds
to traditional smoothed aggregation.
In all experiments involving \spreadlumping, we choose $\tau = 1.1$ as the maximum allowable $\ddom()$ growth.
Results can be found in Table~\ref{tbl:randcube}.  

\begin{table}[htb]\label{tbl:randcube}
\begin{center}
\include{randcube_results_paper}
\end{center}
\caption{
  CG iterations and operator complexities using all combinations of
  enabling and disabling \dmod, \rootstencil, \spreadlumping and \econstraint.
  Each algorithm denoted in \textcolor{red}{red} had 16 solution
  failures due to negative eigenvalues.  These failures are ignored for 
  iteration counts and operator complexities.}
\end{table}

We note first that traditional SA-AMG (all options off) fails in 16 of the 50
test cases.
These failures arise from a negative eigenvalue estimate for $\maxFilEval$ 
in \eqref{eq:prolongator_smoothing} due to the $\filteredA$ matrix, which
has poor diagonal dominance properties. These negative eigenvalues
lead to catastrophic failures within the solver. Here, operator complexity
is defined as the ratio of the number of nonzeros within all the hierarchy
discretization matrices divided by the number of nonzeros for the 
finest level matrix. 
Three algorithms used by themselves (or in combination) --- \dmod,
\rootstencil and \spreadlumping --- prevent all negative eigenvalues in
the test problem.  In the case of \dmod or \spreadlumping, the
diagonal entries are less sensitive to the dropping schemes that lead
to small values in $\filteredD$. The \rootstencil algorithm, 
indirectly removes some problematic entries from $\filteredA$ which
would otherwise lead to negative eigenvalues.  The \econstraint algorithm has no
effect on reducing the number of failures.

Second we note that \spreadlumping and \dmod each modestly reduce the
iteration counts whenever they are used (either individually or in
combination), while \rootstencil tends to slightly increase the
iteration counts on average.  The \econstraint algorithm yields a
substantial decrease in iterations, though it needs to be used in
combination with at least one other option to avoid failures.
For example, the combination of \econstraint and \spreadlumping reduced the iteration
count from $26.1$ with $32\%$ failures to $16.0$ with no failures. This corresponds to
an average iteration reduction of $60\%$ (not counting the failures).
Finally, we note \rootstencil consistently leads to a small reduction
in operator complexity, 1.28 to 1.25, whenever it is used.

\subsection{Triaxially Stretched Cube}\label{sec:results_stretchcube}
We take $(k_x, k_y, k_z)\in \{1,5,10\}^3$, restricted to $k_x
\leq k_y \leq k_z$, yielding 10 different configurations.  For each
$(k_x,k_y,k_z)$ we define a domain $\Omega = [0,3k_x+3] \times [0,3k_y+3] \times [0,3k_z+3]$.  
On $\Omega$, we solve the reaction-diffusion equation
\begin{equation} \label{REACTDIFF_EQN}
\begin{aligned} 
\sigma^{-1} u(x) -\nabla^2 u(x) &= f(x)   \hskip .4in  &x \mbox{~in~} \Omega, \\
u(x) &= u_D(x)                        &x \mbox{~on~} \partial \Omega,
\end{aligned}
\end{equation}
where the forcing function
$f(x)$ is identically zero, 
 $u_D(x) = 1 + x_1 + x_2 + x_3 + x_1 x_2+ x_1 x_3+ x_2 x_3 + x_1 x_2
x_3$
and the reaction term, $\sigma$ is chosen to be $10,10^2,10^3,10^4,$ or $10^5$.
Notice that the nonzero $\sigma$ implies that the matrix has a positive row sum.
We discretize \eqref{REACTDIFF_EQN} using linear hexahedral finite elements on
a tensor product mesh using 60 elements in each coordinate
direction. The element sizes linearly vary from $0.1$ to $k_x/10$ in
the $x$-dimension (and similarly in the $y$ and $z$ dimensions).
Details are given in Appendix~\ref{appendix:stretchcube}.  

Following \S~\ref{sec:results_randcube}, CG with a $10^{-10}$ tolerance is used with a MueLu
preconditioner to solve the linear system.  Again, the distance
Laplacian criterion is used for the aggregation/coarsening phase of
the algorithm with $\theta=0.025$.  Two sweeps of Chebyshev smoothing
is applied (with identical parameters to
\S~\ref{sec:results_randcube}) to all levels except the coarsest
level, where a direct solver is used.

\begin{figure}[htb]
\begin{center}
\subfigure[Traditional SA-AMG (Trad) minus $\dmod$ diagonal modification.]{
    \includegraphics[scale=0.5]{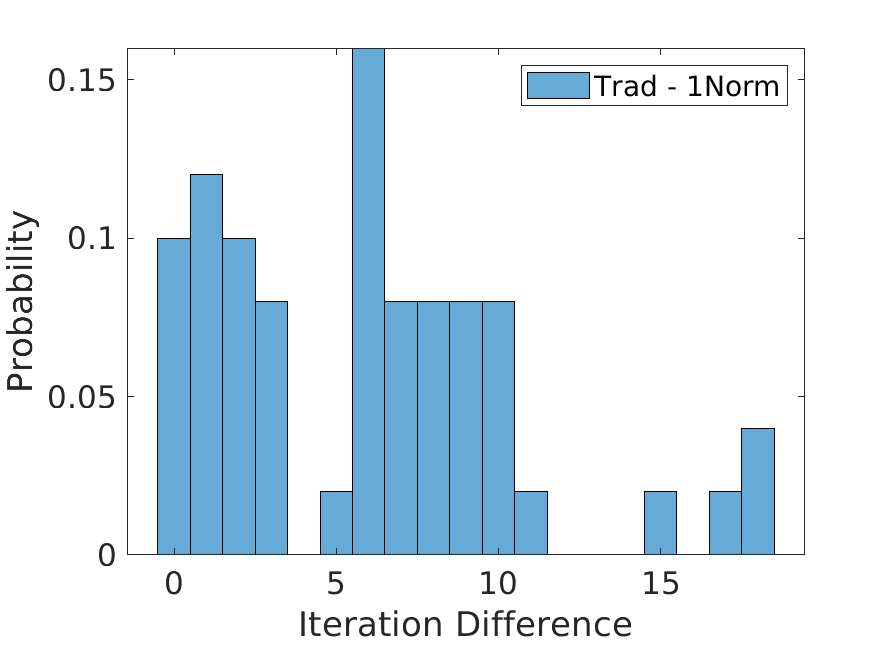}
    \label{fig:stretchcube:dmod}
}
\subfigure[Traditional SA-AMG (Trad) minus $\spreadlumping$ lumping modification.]{
  \includegraphics[scale=0.5]{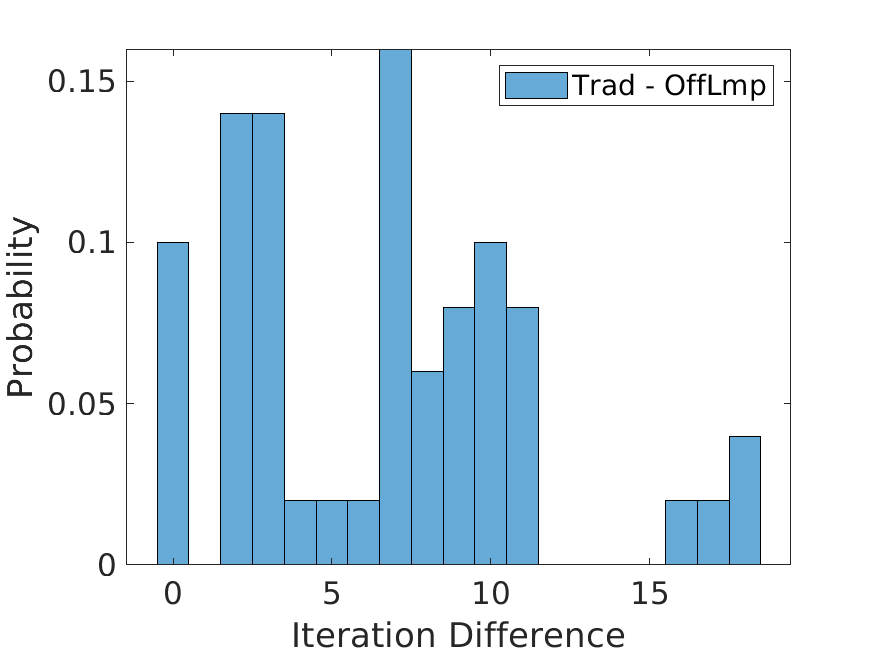}
  \label{fig:stretchcube:spreadlump}
}
\end{center}

\caption{Probability histogram of the difference in iteration counts
  for MueLu-preconditioned CG on \eqref{REACTDIFF_EQN} between
  traditional SA-AMG and one of the proposed algorithms, taken over
  ten mesh configurations $(k_x, k_y, k_z) \in \{1,5,10\}^3$,
  restricted to  $k_x \leq k_y \leq k_z$,
  and five reaction rates $\sigma=1e1,1e2,1e3,1e4,1e5$.}\label{fig:stretchcube}
\end{figure}

First we consider the traditional SA-AMG algorithm and \dmod
modifications. Figure~\ref{fig:stretchcube:dmod} presents a probability
histogram of the number of iterations taken by traditional SA-AMG minus
the number of iterations taken with the \dmod diagonal modification.  We note
that in all fifty cases, the \dmod modification takes no more than the number of
iterations taken by traditional SA-AMG.  While 20 of the \dmod
modification runs are within three iterations of traditional SA-AMG,
four runs have a difference of 15 or more iterations, which is a substantial
savings.  Similarly, Figure~\ref{fig:stretchcube:spreadlump} shows the same
results for the \spreadlumping lumping modification, with 19 of the
runs withing three iterations of traditional SA-AMG, and four runs
with a difference of 15 or more iterations.  We do not show results
for the \econstraint constraint enforcement, but note that 40 of those runs are within three iterations
of traditional SA-AMG while no runs had an iteration difference of 15
or more.  More problematically, two of the 50 runs fail due to
negative eigenvalues. We also do not show results for \rootstencil
sparsification, as in three cases it increases the iteration count substantially.  Thus we
cannot recommend using \econstraint or
\rootstencil on this problem, unless they are paired with
either one of the modifications which can overcome their deficiencies.

\subsection{Subsurface SPE10 Problem}\label{sec:results_sp10}
\begin{figure}
\centering
\includegraphics[scale=0.5]{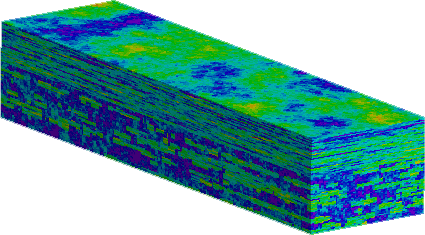}
\caption{Porosity of SPE10 benchmark problem \label{fig:porosity}}
\end{figure}

This model comes from a dataset associated with the $10^{th}$ SPE Comparative Solution Project (SPE10)~\cite{spe10}. It describes a reservoir
simulation that models flow through porous media to predict well production from hydrocarbon deposits. 
The discrete matrix problem was formed using the open source toolbox MRST~\cite{lie_2019}.  The underlying PDE
equations are defined by Darcy's law for a single fluid along with external influences such as wells
and are discretized using a two-point flux approximation method. 
The porosity is shown in figure~\ref{fig:porosity} and leads to large permeability variations that range up to 12 orders of 
magnitude.
The large permeability variation requires that the multigrid solver coarsen irregularly.
For these experiments, the same Chebyshev smoother choices are used as with the first example.
A direct solver is also used on the coarsest grid. A standard smoothed aggregation strength-of-connection criteria is 
employed using the matrix coefficients (i.e., not the distance Laplacian) and three values of 
$\theta$ are considered.  The number of multigrid levels is fixed at 5 where the finest level discretization
matrix is $660K \times 660K $ and the coarsest level matrix dimensions is always less than $14K$
for the largest $\theta$ value and is much smaller for the other two $\theta$.
The conjugate gradient
iteration is terminated when the residual is reduced by a factor of $10^{-8}$.

Table~\ref{results:spe10} illustrates the results for different combinations of our proposed
variants. All combinations that did not use the \dmod variant failed
and so these are not shown.
\begin{table}
\begin{center}
\begin{tabular} {lll||cc|cc|cc}
\toprule
& & & \multicolumn{2}{c|}{$\theta = .02$} &\multicolumn{2}{c|}{$\theta = .05$} &\multicolumn{2}{c}{$\theta = .1$} \\
\multicolumn{3}{c||}{Algorithm Choice}       & AMG      & its.  & AMG      & its.  & AMG      & its. \\
     &   &                                   & complex. &      & complex. &      & complex. &   \\  
\midrule
               &              &              &    1.85  &   24 &  2.00    &   26 &  2.60    &   32 \\ 
               &              & \econstraint &    1.85  &   24 &  2.00    &   26 &  2.60    &   28 \\ 
\spreadlumping &              &              &    1.85  &   23 &  2.00    &   20 &  2.60    &   15 \\ 
\spreadlumping &              & \econstraint &    1.85  &   24 &  2.00    &   20 &  2.60    &   15 \\ 
               & \rootstencil &              &    1.77  &   53 &  1.85    &   35 &  2.12    &   26 \\ 
               & \rootstencil & \econstraint &    1.77  &   49 &  1.85    &   35 &  2.12    &   26 \\ 
\spreadlumping & \rootstencil &              &    1.77  &   48 &  1.85    &   33 &  2.11    &   27 \\ 
\spreadlumping & \rootstencil & \econstraint &    1.77  &   47 &  1.85    &   35 &  2.12    &   27 \\ 
\bottomrule
\end{tabular}
\end{center}
\caption{SPE10 results for combinations of different AMG variants (all results use \dmod).\label{results:spe10}}
\end{table}
For the smallest value of $\theta$, the \econstraint and \spreadlumping procedures do not have too significant
an effect on convergence. The \rootstencil approach does reduce the operator complexity, but
this comes at a fairly significant increase in iteration count (approximately double).  For large values of $\theta$,
however, we see that some of the different algorithm choices do have a more pronounced effect. In particular, the 
best iteration counts employ \spreadlumping for $\theta=.1$, which are about half those of the other methods
when \rootstencil is not used. 
The \rootstencil algorithm more significantly improves the AMG operator complexity when $\theta=.1$ and also improves
the iteration counts in the case that \spreadlumping is not used, though the iterations are worse when both
\spreadlumping and \rootstencil are employed.  The \econstraint option has only a modest effect, giving some 
reduction in iterations for a couple of cases (e.g. when the only other employed variant was \dmod for $\theta=.1$).
The main point, however, is that the traditional SA-AMG algorithm failed to converge on this problem
and so the diagonal $\dmod$ modification is critical in getting the solver to converge.
The additional algorithm modifications can then provide some additional benefit, though not always. 
As the choice of $\theta$ is unknown, it is reassuring that with \spreadlumping the number of 
iterations decrease as the cost/AMG operator complexity increases (associated with an increasing value of $\theta$).

\subsection{Wind Turbine}\label{sec:exawind}

\begin{figure}
\includegraphics[scale=0.15]{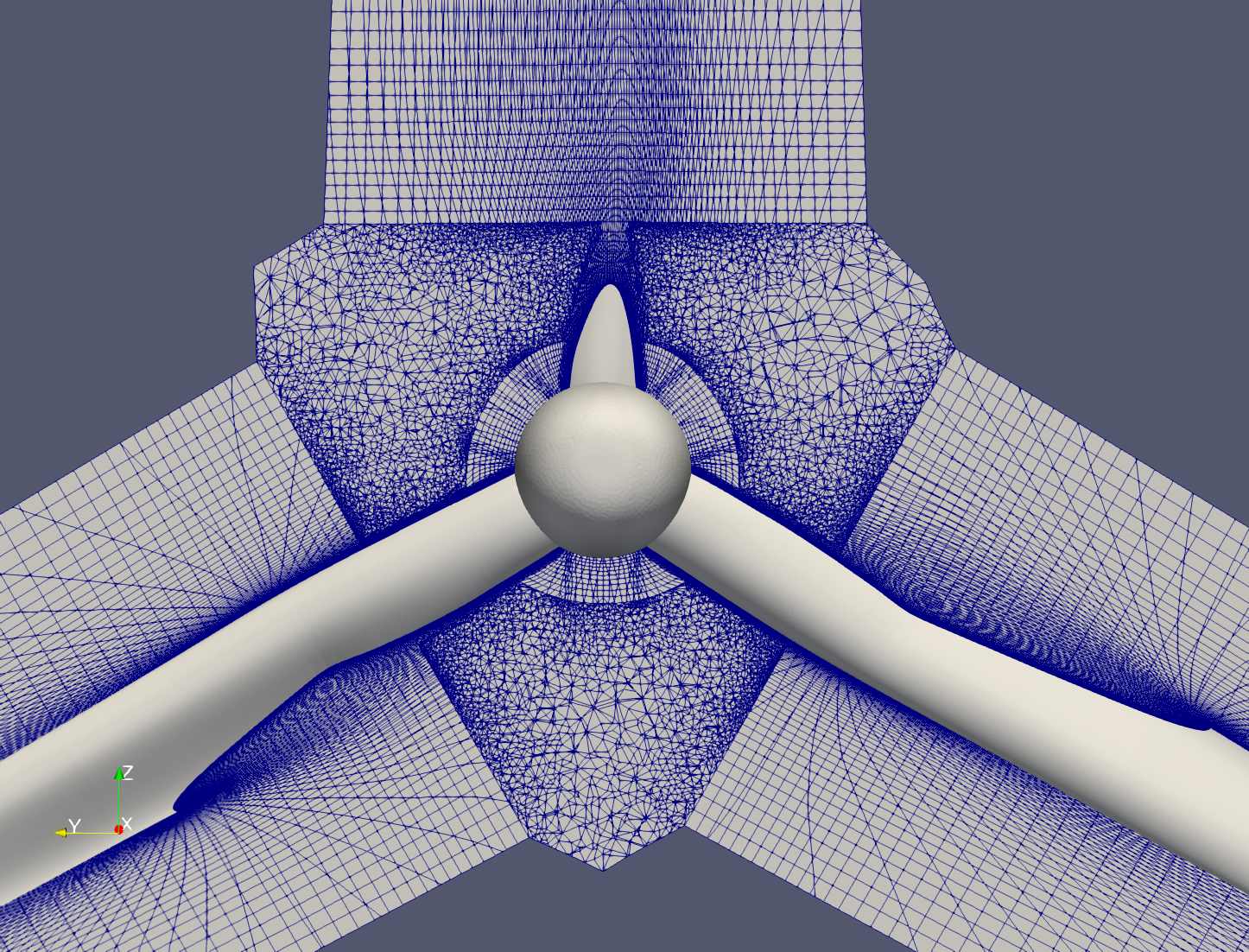}
\includegraphics[scale=0.15]{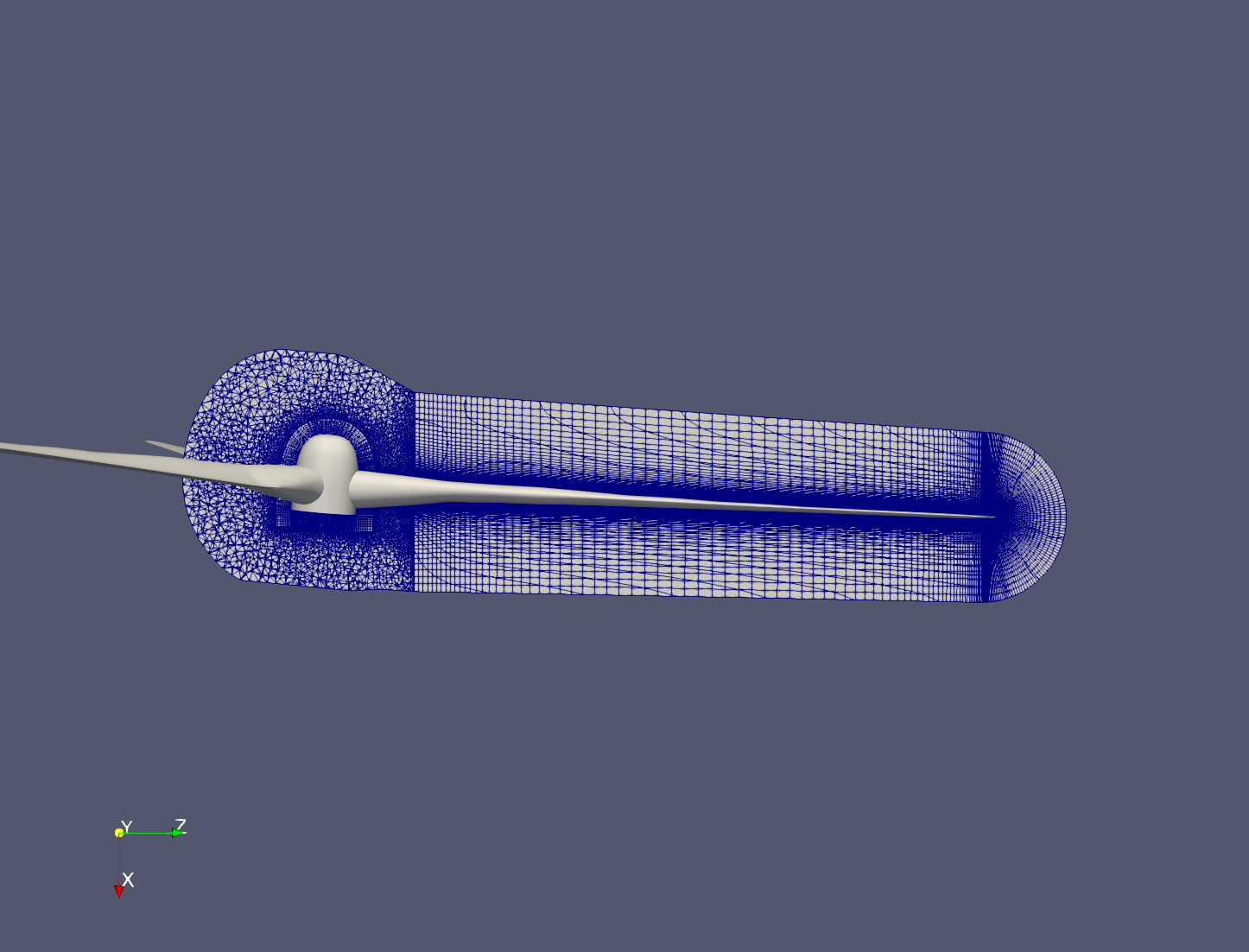}
\caption{Front and side views of NREL 5 megawatt turbine mesh, showing blades and hub.}
\label{fig:nrel5mw}
\end{figure}

We now consider the effects of the various solver options within a wind turbine simulation run in the low Mach computational
fluid dynamics (CFD) code NaluWind~\cite{naluwind}.  Figure~\ref{fig:nrel5mw} depicts a 5-megawatt wind turbine
that is composed of three turbine blades and a hub. 
The blades and background are meshed separately and coupled via constraints.  NaluWind has two main techniques for addressing these constraints.
The first ``coupled" technique can be viewed as an alternating Schwarz approach, in which each mesh has equally valid solutions.
The corresponding linear systems have equations corresponding to constraints, and it is these equations which cause difficulties for the linear solver.
In the second ``decoupled" approach, the constraints are eliminated from the corresponding linear systems, which can be solved independently.

\begin{figure}[tbph]
  \centering
  \includegraphics[scale=0.163]{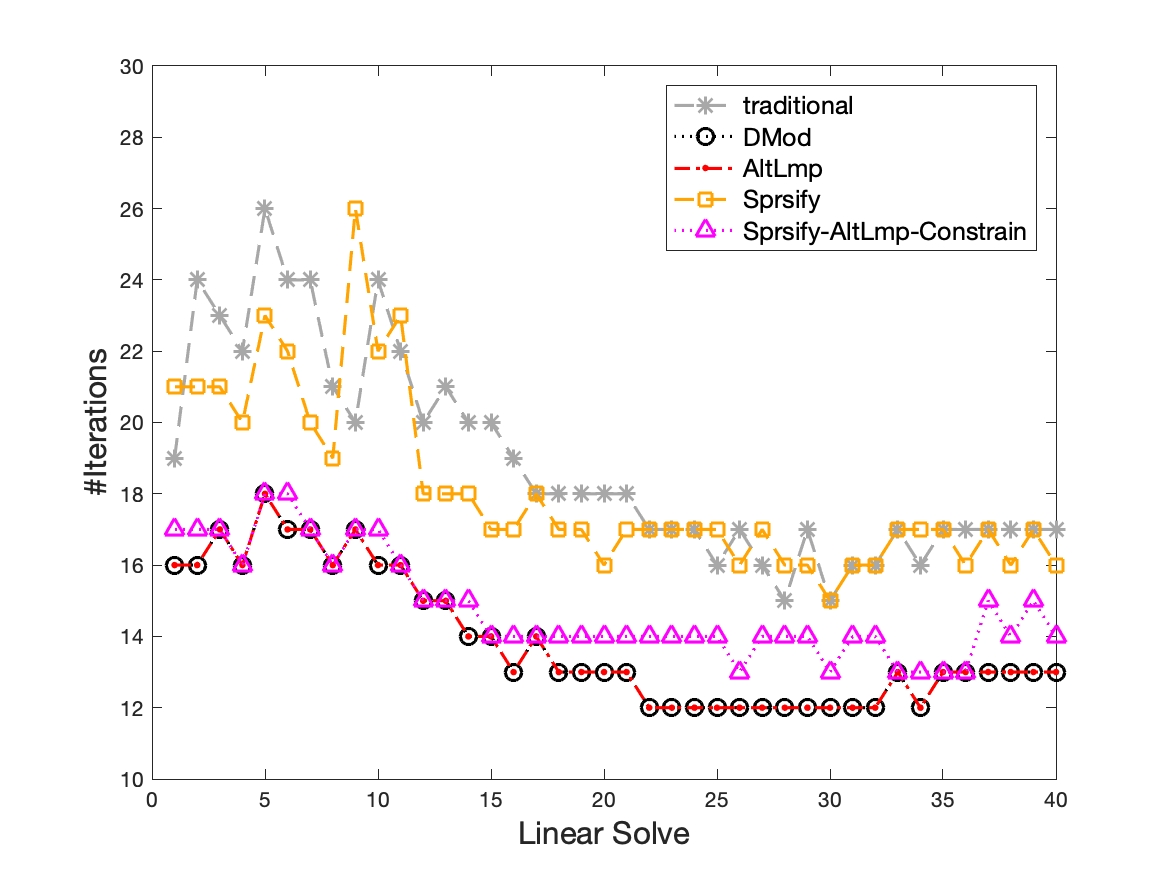}
~~~~~~ \includegraphics[scale=0.45]{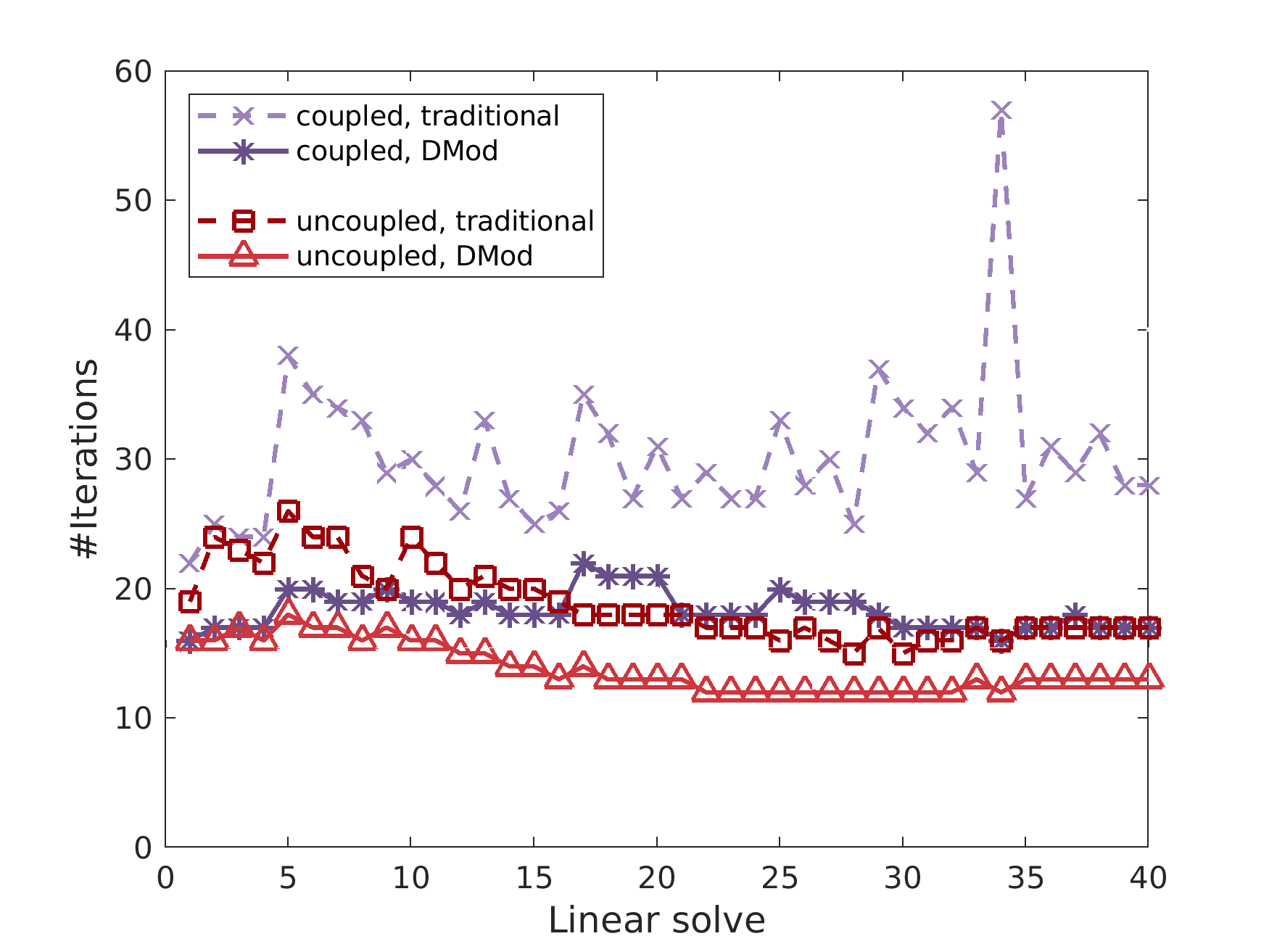} 
  \label{fig:nrel5mw-absrowsum}
  \caption{Comparison of effect of prolongator options on linear solver iterations over 10 timesteps (4 linear solvers per timestep).
           In the left plot, 
           all combinations of the four algorithm variants were tested, but only the most effective strategies are shown.
           ``Traditional" denotes standard SA-AMG.
           ``DMod" is diagonal modification, ``AltLmp" is alternative lumping, and``Sprsify" is prolongator
           sparsification.
           In the right plot, only the effect of "DMod" is shown versus standard SA-AMG for an ExaWind 5MW linear
           system formulated either either as a coupled system (containing constraints) or an uncoupled system.
          }
\end{figure}

\begin{figure}[h!tb]
  \includegraphics[height=0.2\textwidth,width=2.7in]{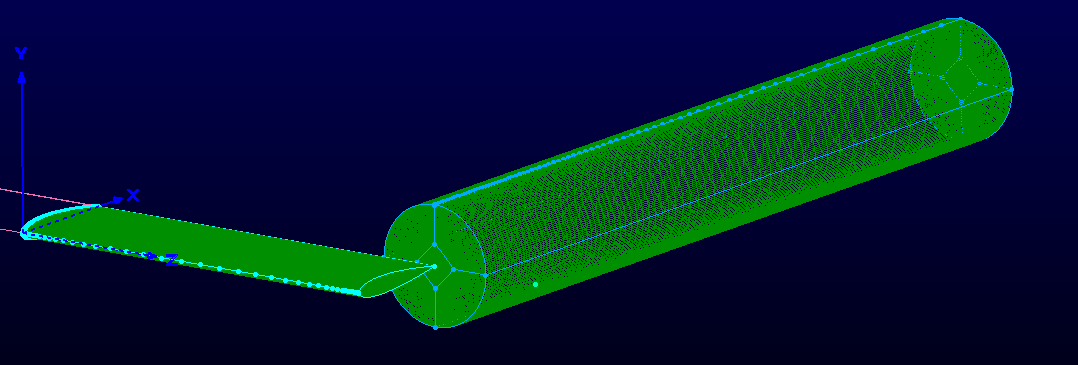}\\
\vskip -1.6in
{\color{white} .}
~~~~~~~~~~~~~~~~~~~~~~~~~~~~~~~~~~~~~~~~~~~~~~~~~~~~~~~~~~~~~~~~~~~~~~
\includegraphics[height=0.27\textwidth]{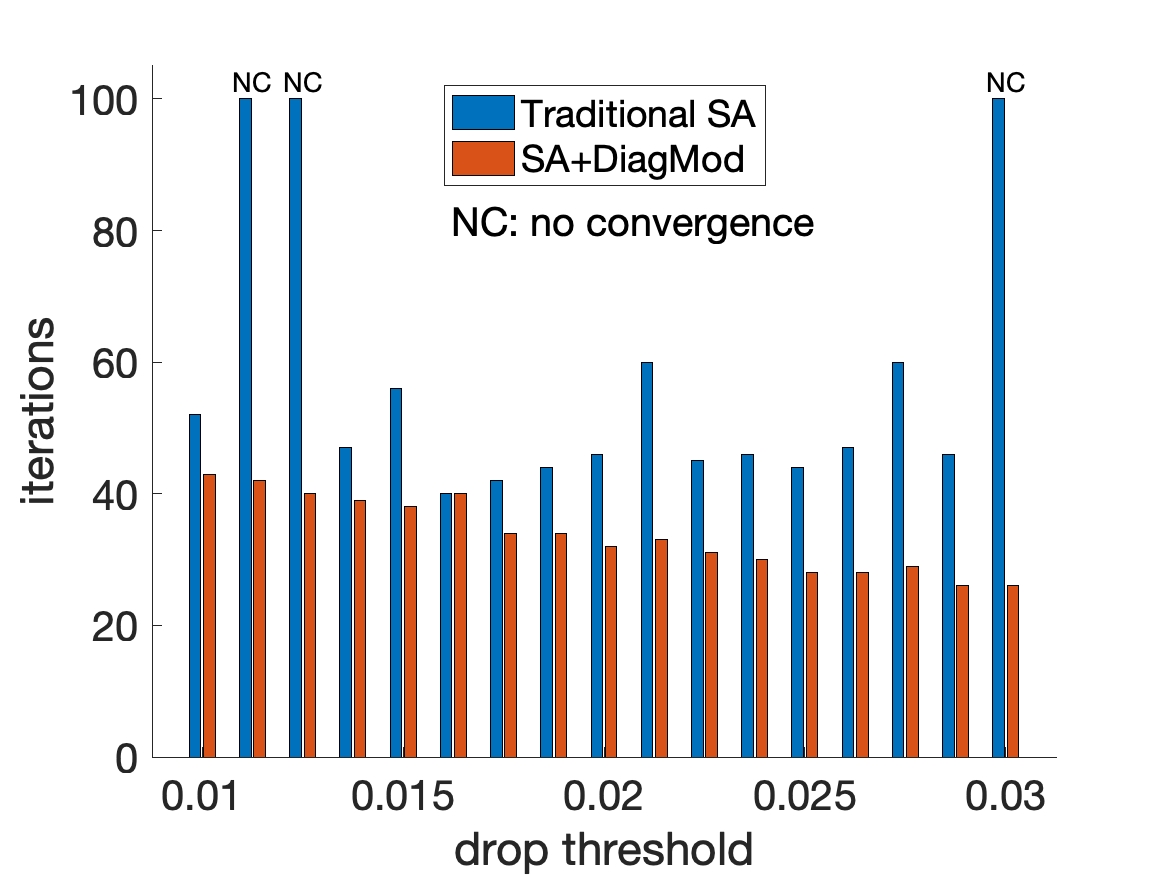}
  \caption{Left: Sample mesh sections for \textsc{5MegaWatt} turbine. Middle: GMRES/AMG iterations for time sequence on 5MegaWatt
   simulation (coupled formulation being more challenging).
   Right: GMRES/AMG iterations for user-supplied strong/weak threshold for a \textsc{McAlister} fixed wing configuration.}
  \label{fig:exawind}
\end{figure}

The simulation itself is time-dependent and consists of two main physics solve phases: momentum and pressure.
While the momentum linear system is amenable to GMRES preconditioned with symmetric Gauss-Seidel, the pressure system
requires a more robust, scalable solver.  The current solver of choice is
GMRES preconditioned by smoothed-aggregation multigrid, with a convergence criteria requiring a relative residual reduction of $10^{-5}$.
For this problem 5 AMG levels are used. The finest level matrix system has 23 million DOFs.   The number of DOFs in
the coarsest level AMG matrices varies, but is typically 9k--15k unknowns. One pre- and one post-smoothing sweep using a degree 2 Chebyshev smoother is employed on all
levels with the exception of the coarsest level where a direct solver is used.  
A standard SA-AMG strength-of-connection criteria is employed using the matrix coefficients (i.e., not the distance Laplacian) for 
a fixed $\theta$ threshold of $.02$. 
We will consider the effect of the proposed SA-AMG options on iteration counts over 10 time steps in the simulation.
The left side of 
Figure~\ref{fig:nrel5mw-absrowsum} compares linear iteration totals for standard SA-AMG versus SA-AMG using various options
for NaluWind run in decoupled mode. The best improvement in iteration counts comes with either the \dmod or the \spreadlumping
options. The \econstraint option is only shown in one case as it generally produced very little improvement.
The \dmod option could not be improved further by using it in conjunction with any of the other algorithms, though it is worth
noting that \spreadlumping and \econstraint did not further degrade performance. As expected, 
\rootstencil generally increased the iteration count and unfortunately provides only a modest gain in AMG operator
complexity from about 1.66 to about 1.61.  When used with the traditional method by itself (shown in the plot) or in 
conjunction with \dmod (not shown), it generally increases the average iteration count by about 1 iteration. 
Finally, when not using the \dmod option, it is possible to improve the \rootstencil iteration counts and make it
fairly competitive with the best runs by applying both \spreadlumping and \econstraint. 
Most of the improvements are due to the \spreadlumping (on average a 2.575 iteration improvement) while the 
\econstraint gains are more modest (on average an additional .75 iteration improvement). 
As the \dmod option alone provides the biggest improvement for this problem, 
the right side highlights the impact of only the \dmod option (all others are turned off) for the coupled formulation. For the sake of 
comparison, the same information is repeated for the decoupled formulation. Here, one can see that the coupled mode leads to generally harder
linear systems. For both coupled and decoupled, there is a nice reduction in the number of iterations using the \dmod option. Further,
the iteration count is generally less erratic as well with the  \dmod option.  Figure~\ref{fig:exawind} considers a robustness
study for one representative linear system within a McAlister fixed wing simulation using the coupled formulation. In this study, 
iteration counts are shown as a function of a varying drop threshold $\theta$ parameter. Without the \dmod option, we again see 
erratic behavior including 3 threshold choices where
the solver does not converge. However, the iteration counts vary smoothly without any failures when the \dmod option is used.
That is, one benefit of \dmod appears to be enhanced robustness.

%% file: randcube_results_paper.tex
\begin{tabular}{|l||lrr|lrr|}
\hline
 &\multicolumn{3}{c|}{\rootstencil off} &\multicolumn{3}{c|}{\rootstencil on}\\
 & &\spreadlumping off &\spreadlumping on & &\spreadlumping off &\spreadlumping on\\
\hline
\multirow{2}{*}{\econstraint off} & \dmod off &\textcolor{red}{26.1(1.28)} &24.2(1.28) & \dmod off &27.9(1.25) &26.1(1.25)\\
 & \dmod on &23.2(1.28) &22.4(1.28) & \dmod on &24.5(1.25) &23.4(1.25)\\
\hline
\multirow{2}{*}{\econstraint on} & \dmod off &\textcolor{red}{19.2(1.28)} &16.0(1.28) & \dmod off &21.2(1.25) &18.2(1.25)\\
 & \dmod on &16.1(1.28) &16.0(1.28) & \dmod on &18.2(1.25) &18.1(1.25)\\
\hline
\end{tabular}

%% file: conclusion.tex
\section{Conclusion} \label{sec:conclusion}
In this paper we have presented four new algorithmic variants to SA-AMG that focus on improving the smoothed prolongator
grid transfer, especially for
problems with many weak connections.  Such systems commonly arise in practice, and can lead to poor SA-AMG performance,
which can manifest as high operator complexity, increased iteration counts, and even failure to converge (due to iteration
matrices with negative eigenvalues).
Whereas other SA-AMG research has focused on developing new strength-of-connection measures to mitigate these issues, we
have assumed a standard scalar strength measure is utilized, and in this paper present algorithms aimed at improving the
final smoothed prolongator.  These algorithms are algebraic in nature and build naturally on the existing SA-AMG machinery.

We have demonstrated the efficacy of these new algorithms on a suite of problems that are challenging for standard SA-AMG to solve:
manufactured scalar Poisson problems with severe variable mesh stretching and lack of diagonal dominance, a standard oil reservoir benchmark, and
linear systems arising from a low-Mach CFD application.
The main take-away is that the four new variants generally yield improvements over standard SA-AMG. The \dmod and \spreadlumping
variations very rarely take more iterations that traditional SA-AMG. While there are cases where the convergence behavior
is similar to that of SA-AMG, there are other cases where \dmod and \spreadlumping are significantly faster and more robust
than SA-AMG. The results with \econstraint are a bit mixed. Sometimes it helps dramatically but other times it is not so
robust. The \rootstencil results do help a modest amount with the multigrid operator complexity, but in most cases
convergence does suffer. However, we do note that \rootstencil was robust on the random cube problem, which was not true 
for SA-AMG. 
Exploring the use of these algorithms in the context of a new weak-connection threshold approach will be the subject of a forthcoming
paper, where \rootstencil's ability to reduce multigrid operator complexity is more significant. 
Another potential topic for future research is the adaptation of one or more of these algorithms to systems of PDEs.

%% file: appendix_randcube.tex
\section{Rand Cube Input}\label{appendix:randcube}
The results from \S~\ref{sec:results_randcube}, were generated from
meshes using the following Pamgen \cite{Pamgen} template.  The mesh is
uniform in $x$ and $y$ and stretched with $100:1$ in the
$z$-direction.  The nodes are each given a random perturbation up to
$20\%$ of the distance to the neighboring node in each direction.
The random number generator was seeded with 50 different seeds in order to
generate the meshes considered.

\begin{verbatim}
mesh
  brick
   zmin =  0.0
   xmin =  0.0
   ymin =  0.0
   numz 1
     zblock 1 1.0 interval 60
   numx 1
     xblock 1 1.0 interval 60
   numy 1
     yblock 1 1.0 interval 60
  end
  set assign
     sideset, ilo, 1
     sideset, jlo, 2
     sideset, klo, 3
     sideset, ihi, 4
     sideset, jhi, 5
     sideset, khi, 6
  end
 user defined geometry transformation
 '
   outxcoord = (inxcoord + 0.2*drand()/60)*1.;
   outycoord = (inycoord + 0.2*drand()/60)*1.0;
   outzcoord = (inzcoord + 0.2*drand()/60)*100.0;
 '
 end
end
\end{verbatim}

%% file: appendix_stretchcube.tex
\section{Triaxially Stretched Cube Input}\label{appendix:stretchcube}
The results from \S~\ref{sec:results_stretchcube}, were generated from
meshes using the following Pamgen \cite{Pamgen} template.  The
elements linear vary in size in each dimension, depending on the
parameters \textsf{KX}, \textsf{KY}, and \textsf{KZ}, which must be
substituted into the input deck below (quantities inside braces are replaced).
\begin{verbatim}
mesh
  brick
   zmin =  0.0
   xmin =  0.0
   ymin =  0.0
   numx 1
     xblock 1 {3.0*(_KX_+1)}, first size .1, last size {_KX_/10}
   numy 1
     yblock 1 {3.0*(_KY_+1)}, first size .1, last size {_KY_/10}
   numz 1
     zblock 1 {3.0*(_KZ_+1)}, first size .1, last size {_KZ_/10}       
  end
  set assign
     sideset, ilo, 1
     sideset, jlo, 2
     sideset, klo, 3
     sideset, ihi, 4
     sideset, jhi, 5
     sideset, khi, 6
  end

end
\end{verbatim}